\documentclass[french, oneside]{amsart}
\usepackage[english,frenchb]{babel}
\usepackage[utf8]{inputenc}
\usepackage[T1]{fontenc}
\usepackage{lmodern}
\usepackage[a4paper]{geometry}
\usepackage{amssymb, amsmath, amsthm}
\usepackage{url}
\usepackage{varioref}
\usepackage[colorlinks, linktocpage, citecolor = red, linkcolor = red, urlcolor=red]{hyperref}
\usepackage{graphicx}
\usepackage{tikz-cd}
\usepackage[shortlabels]{enumitem}
\newtheorem{thm}{Théorème}
\newtheorem{prop}[thm]{Proposition}
\newtheorem{lemme}[thm]{Lemme}

\theoremstyle{definition}

\theoremstyle{remark}
\newtheorem{rmq}[thm]{Remarque}
\newtheorem{ex}[thm]{Exemple}

\newcommand{\sgrn}[1]{\mathopen{\ll} #1\mathopen{\gg}}
\newcommand{\sgr}[1]{\mathopen{<} #1\mathopen{>}}
\newcommand{\N}{\mathbb{N}}

\newcommand{\C}{\mathbb{C}}
\newcommand{\R}{\mathbb{R}}
\renewcommand{\P}{\mathbb{P}^2}
\DeclareMathOperator{\Bir}{Bir}
\DeclareMathOperator{\PGL}{PGL}
\DeclareMathOperator{\Aut}{Aut}
\DeclareMathOperator{\diam}{Diam}

\renewcommand{\H}{\mathbb{H}^{\infty}}

\DeclareMathOperator{\K}{k}

\DeclareMathOperator{\dist}{d}
\newcommand{\longueur}[2]{\dist(#1,#2)}
\newcommand{\isome}[3]{#1_{#2{\scriptscriptstyle{\#}}}^{#3}}
\DeclareMathOperator{\Fix}{Fix}
\DeclareMathOperator{\Tube}{T}
\DeclareMathOperator{\Axe}{Axe}
\DeclareMathOperator{\rayon}{r}
\DeclareMathOperator{\proj}{pr}

\DeclareMathOperator{\PM}{\mathcal{Z}}
\DeclareMathOperator{\NS}{N^1}
\DeclareMathOperator{\cat}{CAT}
\DeclareMathOperator{\Bs}{Bs}
\renewcommand{\epsilon}{\varepsilon}
\newcommand{\plan}{\mathcal{P}}
\newcommand{\A}{\mathbb{A}}
\DeclareMathOperator{\Vois}{Tub}

\DeclareMathOperator{\argcosh}{argcosh}

\title{Non simplicité du groupe de Cremona sur tout corps}
\author{Anne Lonjou}
\address{Institut de Math\'ematiques de Toulouse, Universit\'e Paul Sabatier, 118 route de Narbonne, 31062 Toulouse Cedex 9, France}
\email{alonjou@math.univ-toulouse.fr}



\begin{document}

\begin{abstract}
En utilisant un théorème de F. Dahmani, V. Guirardel et D. Osin, nous prouvons que le groupe de Cremona de dimension $2$ n'est pas simple sur tout corps. Plus précisément, nous montrons que certains éléments de ce groupe satisfont une version affaiblie de la propriété WPD qui est équivalente dans notre contexte particulier à celle de M. Bestvina et K. Fujiwara.
\end{abstract}

\renewcommand{\abstractname}{Abstract}
\begin{abstract}
	Using a theorem of F. Dahmani, V. Guirardel and D. Osin we prove
	that the Cremona group in 2 dimension is not simple,
	over any field. More precisely, we show that some elements
	of this group satisfy a weakened WPD property which is equivalent
	in our particular context to the M. Bestvina and K. Fujiwara's one. (An English version of this paper is available on the webpage of the author.)
\end{abstract}

\maketitle

\section*{Introduction}
Une question importante lors de l'étude d'un groupe est de savoir s'il est simple ou non, et dans ce dernier cas de construire des sous-groupes distingués. Ces questions se sont posées dès la fin du 19\up{ème} siècle pour le groupe $\Bir(\P_{\K})$, le groupe des transformations birationnelles du plan projectif sur un corps $\K$. Cependant, il a fallu attendre 2013 pour que S. Cantat et S. Lamy \cite{CL} y répondent dans le cas où $\K$ est un corps algébriquement clos. \\

L'espace de Picard-Manin associé à $\P_{\K}$ est la limite inductive des groupes de Picard des surfaces obtenues en éclatant toutes suites finies de points de $\P_{\K}$, infiniment proches ou non (voir \cite{Ma}, \cite{C} et le paragraphe \S\ref{subsection_picard_manin}). Il est muni d'une forme d'intersection de signature $(1,\infty)$. En considérant une nappe d'hyperboloïde, nous pouvons lui associer un espace hyperbolique de dimension infinie, noté $\H_{\K}$.
La stratégie de S. Cantat et S. Lamy est de faire agir, par isométries, le groupe de Cremona sur cet espace hyperbolique. Ils obtiennent ainsi de nombreux sous-groupes distingués.
Leur article est constitué de deux parties. 

Dans la première partie, ils définissent la notion d'élément tendu (\og tight element\fg \  dans leur terminologie) que nous énonçons dans le cas particulier du groupe de Cremona, de la manière suivante : un élément $g$ de $\Bir(\P_{\K})$ est tendu si l'isométrie correspondant à $g$ est hyperbolique et vérifie deux conditions. La première est que son axe soit rigide, c'est-à-dire si l'axe d'un conjugué de $g$ est proche de l'axe de $g$ sur une distance suffisamment grande alors les deux axes sont confondus. La seconde condition est que le stabilisateur de l'axe de $g$ soit le normalisateur du groupe $\sgr{g}$. 
Plus précisément, l'élément hyperbolique $g$ doit vérifier : \begin{enumerate}[label=\textbullet,font=\small,leftmargin=1cm,wide]
	\item Pour tout $\epsilon>0$, il existe $C\geq 0$ tel que si un élément $f$ dans $\Bir(\P_{\K})$ satisfait l'inégalité $\diam(\Vois_{\epsilon}(\Axe(g))\cap \Vois_{\epsilon}(f\Axe(g)))\geq C$ alors $f\Axe(g)=\Axe(g)$, où $\Vois_{\epsilon}$ signifie le voisinage tubulaire de rayon $\epsilon$.
	\item Pour tout élément $f$ dans $\Bir(\P_{\K})$, si $f\Axe(g)=\Axe(g)$ alors $fgf^{-1}=g^{\pm1}$.
\end{enumerate} 
Ce critère leur permet d'établir une variante de la propriété de la petite simplification : 

\begin{thm}\label{theo_CL}\cite{CL}
	Soit $\K$ un corps algébriquement clos. Si $g\in\Bir(\P_{\K})$ est tendu alors il existe un entier $n$ non nul tel que pour tout élément $h$ non trivial appartenant au sous-groupe normal engendré par $g^n$, le degré de $h$ vérifie $\deg(h)\geq \deg(g^n)$. En particulier, le sous-groupe normal $\sgrn{g^n}$ de $\Bir(\P_{\K})$ est propre. 
\end{thm}
La seconde partie consiste à montrer qu'il existe des éléments tendus dans $\Bir(\P_{\K})$. Nous y reviendrons à la fin de cette introduction.\\

L'objectif du présent article est d'obtenir une preuve de la non simplicité du groupe de Cremona qui fonctionne pour tout corps $\K$ : 
\begin{thm}\label{theo_principal}
	Pour tout corps ${\K}$, le groupe de Cremona $\Bir(\P_{\K})$ n'est pas simple.
\end{thm}

Pour démontrer ce théorème, nous n'allons pas utiliser le fait qu'un élément soit tendu mais plutôt le fait qu'il satisfasse la propriété WPD (\og weakly properly discontinuous\fg), propriété proposée par M. Bestvina et K. Fujiwara \cite{BF} en 2002 dans le contexte du mapping class group. Un élément $g$ satisfait la propriété WPD, si pour tout $\epsilon\geq 0$, il existe un point $x$ et un entier $n$ strictement positif tels qu'il n'existe qu'un nombre fini d'éléments du groupe $G$ déplaçant $x$ et $g^n(x)$ d'au plus $\epsilon$. Les éléments hyperboliques que nous étudions possèdent un axe. Par conséquent, nous emploierons de manière équivalente la terminologie introduite par R. Coulon dans son exposé au séminaire Bourbaki \cite{Co} : le groupe $G$ agit discrètement le long de l'axe de $g$. L'avantage de la formulation de R. Coulon est d'expliciter le rôle du groupe $G$.
À noter que D. Osin \cite{Os} a englobé cette notion ainsi que d'autres sous le terme d'action acylindrique. Il unifie ainsi plusieurs travaux concernant des groupes différents. \\

Récemment, F. Dahmani, V. Guirardel et D. Osin \cite{DGO} généralisent également la théorie de la petite simplification pour des groupes agissant par isométries sur des espaces $\delta$-hyperboliques. Rappelons qu'un espace métrique géodésique $X$ est $\delta$-hyperbolique si pour chaque triangle de $X$, tout côté est contenu dans le $\delta$-voisinage de la réunion de ses deux autres côtés. Une de leur motivation est d'étudier le mapping class group d'une surface de Riemann hyperbolique. Ce groupe agit sur le complexe des courbes qui est un espace $\delta$-hyperbolique non localement compact (tout comme l'espace $\H_{\K}$ mentionné plus haut). Dans ce contexte, ils construisent des sous-groupes distingués propres qui sont de plus libres et purement pseudo-Anosov, c'est-à-dire que tous leurs éléments non triviaux sont pseudo-Anosov. Ils répondent ainsi à deux questions restées longtemps ouvertes. 

Le lien entre la théorie de la petite simplification et la propriété WPD se fait au travers de deux énoncés \cite[Theorem 1.3 et Corollary 2.9]{Gui}. Le premier nous dit que dans le groupe normal engendré par une famille vérifiant la propriété de petite simplification, les éléments ont une grande longueur de translation, et le second nous dit que lorsqu'un élément $g$ satisfait la propriété WPD, les conjugués de $\sgr{g^n}$ forment une famille satisfaisant la propriété de petite simplification. En réunissant ces deux énoncés (voir également \cite[Theorem 5.3, Proposition 6.34]{DGO}), nous obtenons le théorème :

\begin{thm}\label{theo_DGO}  
	\cite{DGO}
	Soit $C$ un nombre réel positif. Soient $G$ un groupe agissant par isométries sur un espace $X$ $\delta$-hyperbolique et $g$ un élément hyperbolique de $G$. Si $G$ agit discrètement le long de l'axe de $g$ alors il existe $n\in \N$ tel que pour tout élément $h$ non trivial et appartenant au sous-groupe normal engendré par $g^n$, $L(h) > C$ où $L$ est la longueur de translation. En particulier, pour $n$ assez grand, le sous-groupe normal $\sgrn{g^n}$ de $G$ est propre. De plus, ce sous-groupe est libre.
\end{thm}

En utilisant ce théorème, la preuve du théorème \ref{theo_principal} se résume à exhiber des éléments satisfaisant la propriété WPD. Nous trouvons de tels éléments dans $\Aut(\A_{\K}^2)$ que nous identifions à un sous-groupe de $\Bir(\P_{\K})$ via l'application qui envoie $(x,y)\in\A_{\K}^2$ sur $[x:y:1]\in\P_{\K}$. Notre résultat principal est alors :
\begin{thm}\label{prop_principale}
	Soient $n\geq 2$ et $\K$ un corps de caractéristique ne divisant pas $n$. Considérons l'action du groupe $\Bir(\P_{\K})$ sur $\H_{\bar{{\K}}}$ où $\bar{\K}$ est la clôture algébrique de $\K$.  Le groupe $\Bir(\P_{\K})$ agit discrètement le long de l'axe de l'application : \[\begin{array}[t]{lrcl}
	h_n : & \A_{\K}^2 & \longrightarrow & \A_{\K}^2 \\
	& (x,y) & \longmapsto & (y,y^n-x) \end{array}.\]
\end{thm}

Remarquons que si ${\K}$ est un corps algébriquement clos de caractéristique $p>0$, pour tout entier $k\geq 1$, le sous-groupe normal engendré par $h_p^k$ est le groupe $\Bir(\P_{\K})$. En effet, $h_p$ normalise les translations : 
\[(x^p-y,x)\circ(x+a,y+b)\circ(y,y^p-x)=(x+a^p-b,y+a).\]
Nous obtenons que $\sgrn{h_p^k}$ contient une translation non triviale. Le théorème de Noether permet alors de montrer l'égalité $\sgrn{h_p^k}=\Bir(\P_{\K})$, voir par exemple \cite[Proposition 5.12]{CD} ou \cite[p.42]{Gi}. Plus généralement, sur un corps infini de caractéristique divisant $n$, l'application $h_n$ ne satisfait pas la propriété WPD, ce qui justifie l'énoncé du théorème \ref{prop_principale}.

En fait, comme conséquence des résultats de F. Dahmani, V. Guirardel et D. Osin, nous obtenons des propriétés plus précises que la seule non-simplicité :
\begin{thm}
	Soit $k$ un corps. Le groupe $\Bir(\P_k)$ contient des sous-groupes distingués libres, et est SQ-universel.
\end{thm}
Rappelons qu'un groupe $G$ est SQ-universel si tout groupe dénombrable se plonge dans un quotient de $G$. Plus de détails se trouvent par exemple dans \cite[Theorem 2.14]{Gui}. \\

Nous terminons cette introduction en comparant notre preuve du théorème \ref{theo_principal} à des résultats récents de \cite{SB} d'une part, et à la stratégie de \cite{CL} d'autre part.

Dans son article \cite[Corollary 7.11]{SB}, N.I. Shepherd-Barron prouve que tout élément hyperbolique du groupe de Cremona sur un corps ${\K}$ fini engendre un sous-groupe normal propre. En particulier, $\Bir(\P_{\K})$ n'est pas simple pour tout corps ${\K}$ fini. Dans ce même article (Theorem 7.6), il donne un critère en terme de la longueur de translation de $g$, pour qu'une transformation hyperbolique $g$ soit tendue, et donc que le sous-groupe distingué engendré par l'une de ses puissances soit propre.
\begin{thm}\cite{SB}
Supposons que la caractéristique de $\K$ soit nulle ou que $\K$ soit de caractéristique $p>0$ et algébrique sur $\mathbb{F}_p$.
Si la longueur de translation d'un élément hyperbolique $g$ de $\Bir(\P_{\K})$ n'est pas le logarithme d'un entier quadratique, ni de la forme $\log p^n$ dans le cas où la caractéristique est strictement positive, alors une certaine puissance de $g$ est tendue.
\end{thm}

 Une hypothèse sans doute excessive est faite sur le corps dans le cas de la caractéristique positive (${\K}$ doit être isomorphe à un sous-corps de $\bar{\mathbb{F}_p}$), dans le but d'éviter le problème mentionné ci-dessus des transformations normalisant le sous-groupe des translations. Cependant, même si nous pouvions enlever cette hypothèse et obtenir ainsi une démonstration alternative de la non simplicité de $\Bir(\P)$ sur tout corps, une telle preuve reste peu élémentaire car elle repose lourdement sur les articles \cite{CL} et \cite{BC}.

Dans ce dernier article, J. Blanc et S. Cantat  s'intéressent aux degrés dynamiques de toutes les applications birationnelles de surface projective. Ils prouvent notamment qu'il n'y a pas de degré dynamique dans l'intervalle $\left]1,\lambda_L\right[$ où $\lambda_L$ est le nombre de Lehmer (\og gap property\fg). Ils obtiennent comme corollaire que pour tout élément hyperbolique $g$ du groupe de Cremona, l'indice de $\sgr{g}$ dans son centralisateur est fini. Comme remarqué par R. Coulon \cite{Co}, ceci implique que si un élément est tendu alors $\Bir(\P_{\K})$ agit discrètement le long de l'axe de cet élément. 

Dans \cite{CL} une relation entre l'entier $n$ du théorème \ref{theo_CL} et la longueur de translation de $g$ est donnée. Une autre conséquence de la \og gap property\fg\ est que l'entier $n$ peut être choisi de façon uniforme : $n\geq \max \{\frac{139347}{\lambda_L},\frac{10795}{\lambda_L}+374\}$ convient.\\

Dans leur article \cite{CL}, S. Cantat et S. Lamy exhibent de deux manières différentes des éléments tendus, selon si le corps est $\C$ ou un corps algébriquement clos quelconque. 

Dans le cas où ${\K}$ est un corps algébriquement clos, ils éclatent des points particuliers de $\P_{\K}$ pour obtenir une surface $S$ ayant des automorphismes de grands degrés dynamiques. Par exemple, éclater les $10$ points doubles d'une sextique rationnelle donne une telle surface, dite de Coble. 
Si le corps n'est pas algébriquement clos, les automorphismes sont à coefficients dans la clôture algébrique mais pas dans le corps ${\K}$ initial, leur preuve ne s'étend donc pas dans ce cas là.

Dans le cas où le corps est $\C$, ils considèrent un élément \og général\fg\  du groupe de Cremona de la forme $g=a \circ J$ où $a\in \PGL_3(\C)$ et $J$ est une transformation de Jonquières, et ils montrent qu'il est tendu. Le qualificatif \og général\fg \  signifie que n'importe quel élément $a\in \PGL_3(\C)$ convient, quitte à enlever un nombre dénombrable de fermés de Zariski propres dans $\PGL_3(\C)$. Ainsi si le corps de base est dénombrable, il se peut qu'il n'y ait pas d'application générale. C'est pourtant cette méthode que nous généralisons. En effet, dans leur preuve, il est nécessaire que $g$ satisfasse deux conditions. La première impose aux points base de $g$ et de $g^{-1}$ d'être dans $\P_{\K}$. La seconde impose aux points base des itérés de $g$ d'être disjoints de ceux des itérés de $g^{-1}$. Les applications $h_n$ que nous considérons se décomposent comme la composée d'une involution de Jonquières et de l'involution linéaire qui échange les coordonnées, $h_n=(y,x)\circ(y^n-x,y)$. La première condition n'est pas vérifiée car $h_n$ et $h_n^{-1}$ n'ont qu'un seul point base dans $\P_{\K}$ mais nous verrons que ce n'est pas un problème. Concernant la seconde condition, composer par l'élément $(y,x)\in\PGL_3({\K})$ permet de séparer les points base des itérés de $h_n$ de ceux des itérés de $h_n^{-1}$. Ainsi cet élément particulier joue le même rôle que l'élément général $a$ ci-dessus. La preuve présentée dans cet article ne repose pas sur \cite{BC,CL} à l'exception du lemme~\ref{lemme_lineaire} à la section \ref{section_resprincipal} qui est une adaptation directe de \cite[Proposition 5.7]{CL}.\\

L'article est organisé comme suit. Dans la première section, nous rappelons la construction formelle d'un espace hyperbolique de dimension infinie puis nous introduisons la notion de \og tube géodésique\fg, qui nous permet notamment de court-circuiter \cite[Lemma 3.1]{CL}. Le résultat principal de cette section permet d'affaiblir, dans un espace hyperbolique de dimension infinie, les hypothèses assurant à un groupe d'agir discrètement le long de l'axe d'un de ses éléments. Dans la seconde section, nous rappelons la définition de l'espace de Picard-Manin ainsi que la construction de l'espace hyperbolique $\H_{\bar{\K}}$ associé. Nous utilisons l'action du groupe $\Bir(\P_{\K})$ sur $\H_{\bar{\K}}$ afin d'exhiber des éléments satisfaisant la propriété WPD. Ceci nous permet de prouver le théorème \ref{prop_principale} et par conséquent de montrer que le groupe de Cremona n'est pas simple pour tout corps.

\section*{Remerciements}
Je remercie vivement Stéphane Lamy, mon directeur de thèse,
qui m'a proposé ce sujet. Je le remercie également pour ses nombreuses relectures.	
Merci au rapporteur pour ses commentaires qui ont permis de rendre certains passages plus clairs.

\section{Propriété WPD dans un espace hyperbolique de dimension infinie}\label{section_hyperbolique}

L'objectif de cette section est de montrer que lorsqu'un groupe agit sur un espace hyperbolique de dimension infinie, nous pouvons affaiblir les hypothèses qui assurent qu'un groupe agit discrètement le long de l'axe d'un de ses éléments.	
\subsection{Espace hyperbolique de dimension infinie}
Nous rappelons ici la construction des espaces hyperboliques de dimension infinie (voir par exemple \cite[Section 6.3]{C}).
Soit $(H,\langle \cdot,\cdot\rangle)$ un espace de Hilbert réel de dimension infinie. Soient $u\in H$ de norme $1$ et $u^{\perp}$ son complémentaire orthogonal. Tout élément $v\in H$ s'écrit de manière unique $v=v_uu+v_{u^\perp}$ où $v_u\in\R$ et $v_{u^\perp}$ appartient à $u^{\perp}$. Nous définissons une forme bilinéaire symétrique $\mathcal{B}$ sur $H$ de signature $(1,\infty)$ de la manière  suivante : $\mathcal{B}(x,y)=x_uy_u-\langle x_{u^{\perp}}, y_{u^{\perp}}\rangle$. Notons $\H$ la nappe d'hyperboloïde définie par : \[ \H:=\{x\in H\mid \mathcal{B}(x,x)=1 \text{ et } \mathcal{B}(u,x)>0 \}.\]

L'espace $\H$ muni de la distance $\dist$ définie par $\cosh\dist(x,y) :=\mathcal{B}(x,y)$ est un espace métrique complet et de dimension infinie. Remarquons que si l'intersection de $\H$ avec un sous-espace vectoriel de dimension $n+1$ de $H$ est non vide alors c'est une copie de l'espace hyperbolique usuel $\mathbb{H}^n$. En particulier, il existe une unique géodésique reliant deux points de $\H$, obtenue comme l'intersection de $\H$ avec le plan vectoriel contenant ces deux points. Dans la pratique, nous nous ramenons souvent à $\mathbb{H}^2$ en considérant un plan hyperbolique de $\H$. De ce fait, tout triangle de $\H$ est isométrique à un triangle de $\mathbb{H}^2$. Ceci implique que $\H$ est $\cat(-1)$ et $\delta$-hyperbolique pour la même constante $\delta =\log(1+\sqrt{2})\approx 0.881373587$ que $\mathbb{H}^2$ (voir \cite[Example 1.23]{Cal}).
Dans la suite lorsque nous parlerons d'un espace hyperbolique $\H$, il sera toujours de cette forme.

Introduisons à présent quelques définitions et notations.
Soit $f$ une isométrie de $\H$, sa longueur de translation est définie par $L(f)=\underset{x\in\H}{\inf}\dist(x,f(x))$. Si la longueur de translation de $f$ est strictement positive, $f$ est dite hyperbolique. Dans ce cas, elle possède un axe invariant déterminé par les points réalisant l'infimum : \[\Axe(f):=\{x\in\H\mid \dist(x,f(x))=L(f) \}.\]
De plus, $f$ s'étend de manière unique au bord $\partial\H$. Pour plus de détails sur la notion de bord, nous renvoyons à \cite[p.260]{BH}. L'isométrie hyperbolique $f$ possède exactement deux points fixes sur $\partial\H$ ; l'un répulsif noté $b^-$, et l'autre attractif noté $b^+$. Ce sont les bouts de l'axe de $f$ que l'on oriente de $b^-$ vers $b^+$. Nous définissons ainsi une relation d'ordre sur les points de $\Axe(f)$. Le point $x$ est dit plus petit que $y$, noté $x<y$ si $x\in\left]b^-,y\right]$ (et par symétrie $y\in [x,b^+[$). Lorsqu'un point $x$ est sur l'axe de $f$ nous notons respectivement $x-\epsilon$ et $x+\epsilon$ les deux points sur $\Axe(f)$ situés à distance $\epsilon$ de $x$ tels que $x-\epsilon<x+\epsilon$.
Remarquons que pour tout point $x\in\H$, la suite $(f^{\pm n}(x))_{n\in\N}$ converge vers $b^{\pm}$.

\subsection{Propriété WPD}
	
Soit $G$ un groupe qui agit par isométrie sur un espace métrique $(X,\dist)$. Pour toute partie $A$ de $X$ et pour toute constante $\epsilon\geq 0$, nous notons \[\Fix_{\epsilon}A:=\{g\in G\mid \dist(a,ga)\leq \epsilon \ \ \forall a\in A \},\] le fixateur à $\epsilon$ près de $A$ par $G$. Dans la suite nous ne préciserons plus que les actions considérées sont isométriques.
Le lemme suivant est bien connu.

\begin{lemme}\label{lemme_WPD_x}
	Soient $G$ un groupe agissant sur un espace métrique $X$ et $g$ un élément de $G$. Les deux propriétés suivantes sont équivalentes :\begin{enumerate}
		\item\label{lemme_point_1} Il existe $y\in X$ tel que pour tout $\epsilon\geq 0$, il existe $N\in\N$ tel que $\Fix_{\epsilon}\{y,g^N(y)\}$ soit fini.
		\item Pour tout $x\in X$, pour tout $\epsilon\geq 0$, il existe $N\in\N$ tel que $\Fix_{\epsilon}\{x,g^N(x)\}$ soit fini.
	\end{enumerate}
\end{lemme}
\begin{proof}
	Soit $y$ vérifiant le point (\ref{lemme_point_1}). Soient $x$ un point quelconque de $X$ et $\epsilon\geq 0$. Posons $\epsilon'=2\dist(x,y)+\epsilon$ et $N$ l'entier du point (\ref{lemme_point_1}) associé à $\epsilon'$. Par l'inégalité triangulaire, nous avons l'inclusion :\[\Fix_{\epsilon}\{x,g^N(x)\}\subset\Fix_{\epsilon'}\{y,g^N(y)\}.\] Ainsi, l'ensemble $\Fix_{\epsilon}\{x,g^N(x)\}$ est fini.
\end{proof}

Nous disons dans la situation du lemme que $g$ \textit{satisfait la propriété WPD} (\og weak proper discontinuity\fg), ou encore lorsque la notion d'axe est bien définie que le groupe $G$ \textit{agit discrètement le long de l'axe de $g$} (voir \cite{BF} et \cite{Co}).\\

Dans les prochains paragraphes, nous introduisons la notion de \og tube géodésique \fg. Ceci nous permet d'affaiblir les hypothèses garantissant la propriété WPD dans le contexte d'un groupe agissant sur $\H$.

\subsection{Quadrilatères hyperboliques}
    Nous rappelons ici un lemme donnant des relations trigonométriques pour des quadrilatères hyperboliques ayant trois angles droits.
	
	\begin{lemme}\label{lemme_relation-quadrilatere}
	Soit $ADCB$ un quadrilatère de $\mathbb{H}^2$ dont les angles en $B$, $C$ et $D$ sont droits, alors la relation suivante est vérifiée : \[\tanh\longueur{A}{B}=\tanh\longueur{D}{C}\cosh\longueur{C}{B}.\]  
	\end{lemme}	
	\begin{figure}[h]
		\centering
		\def\svgwidth{150pt}
		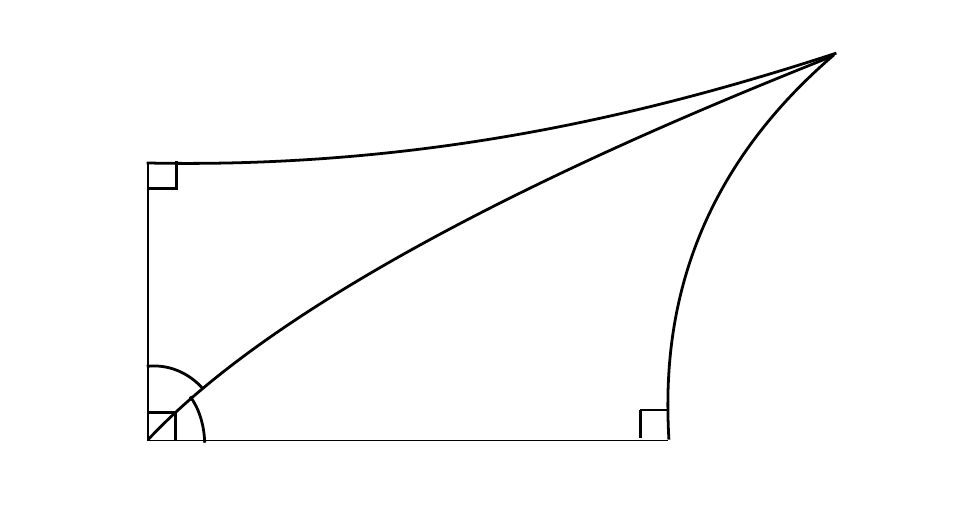
		\caption {Quadrilatère ADCB.\label{figure_quadrilatère}}
	\end{figure}
	
	\begin{proof}
	Dans le triangle rectangle $ABC$, notons $a$, $b$, $c$ les longueurs opposées aux sommets $A$, $B$, $C$, et $\gamma$ l'angle en $C$, et dans le triangle $ACD$, $d$ et $e$, les longueurs opposées aux sommets $A$ et $C$ (voir Figure \ref{figure_quadrilatère}).
	 Les triangles rectangles $ABC$ et $ACD$ satisfont les relations classiques (voir par exemple \cite[Theorem 7.11.2]{B}) : \[ \cos \gamma=\frac{\tanh a}{\tanh b}\ \  \text{,   } \ \sin \gamma=\cos (\frac{\pi}{2}-\gamma)=\frac{\tanh d}{\tanh b} \text{ et } \tan\gamma=\frac{\tanh c}{\sinh a}.\] 

	Cela nous donne l'égalité cherchée : \begin{equation*}
		\tanh c=\tanh d\cosh a. \qedhere
	\end{equation*}
	\end{proof}
	\subsection{Tubes géodésiques}
	Fixons un espace $(\H, \dist)$ comme défini précédemment. Soit $\Gamma$ une géodésique de $\H$. Nous introduisons à présent la notion de \textit{tube géodésique} autour de $\Gamma$ qui sera centrale dans la preuve de la proposition principale de cette section (Proposition \ref{pro_equiv_WPD}). Le mot \og géodésique \fg\  est employé ici pour signifier que lorsque l'on considère l'intersection de ce tube avec un plan hyperbolique $\plan$ contenant $\Gamma$, le quadrilatère plein obtenu est bordé de segments géodésiques.

\begin{figure}
		\hspace{-5cm}
		\def\svgwidth{330pt}
		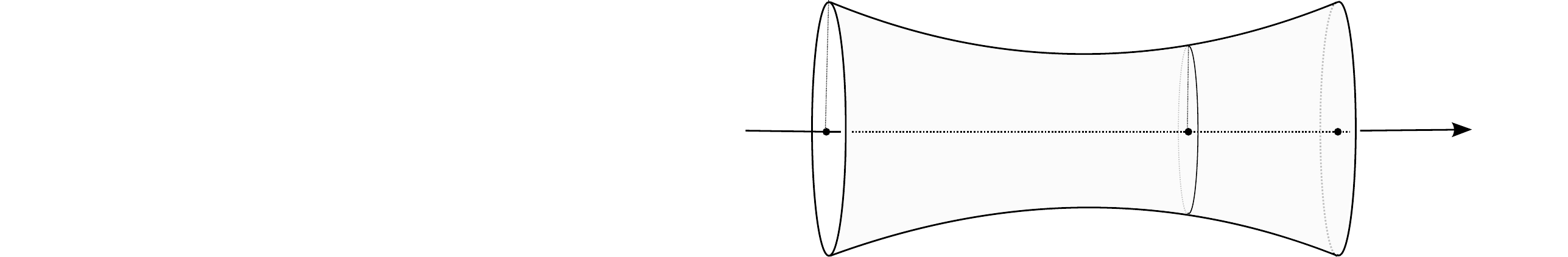
		\caption {Tube $\Tube_{x,x'}^{\epsilon}$. \label{figure_tube}}
		\end{figure}

\begin{figure}
		 				\centering
		 				\def\svgwidth{290pt}
		 				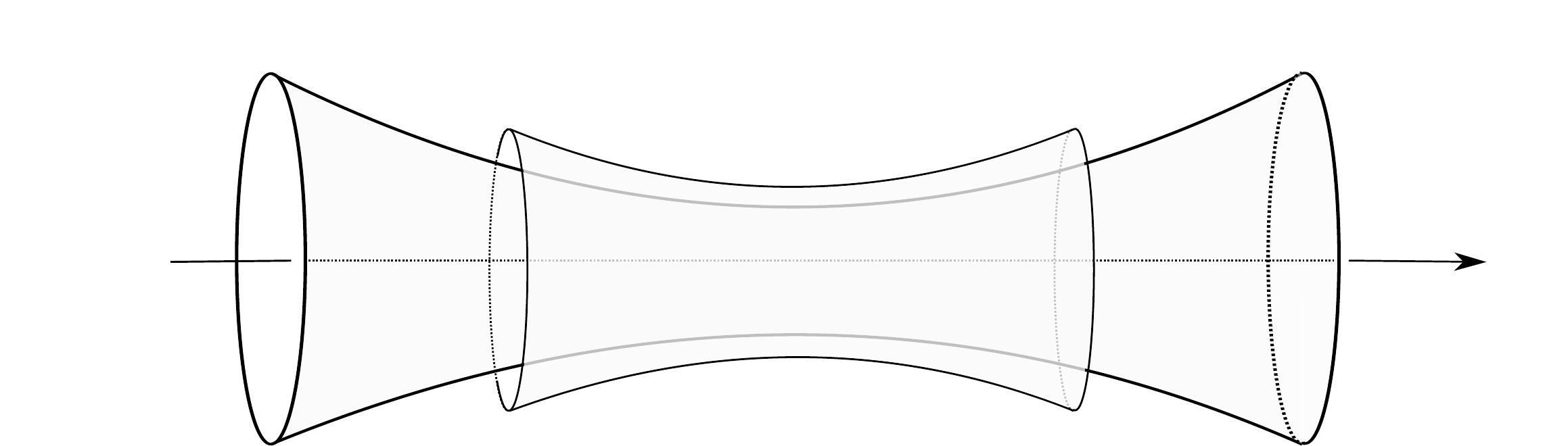
		 			\caption {Tube $\Tube_{x,x'}^{\epsilon_1}$ traversant le tube $\Tube_{y,y'}^{\epsilon_2}$.\label{figure_tube_passant}}
		 		\end{figure}

	Soient $x$ et $x'$ deux points de $\H$ et $\Gamma$ la géodésique passant par ces deux points. Soit $\proj_{\Gamma}$ la projection qui envoie un point $y\in\H$ sur le point de $\Gamma$ le plus proche de $y$. Nous notons \[\Gamma^{\perp}_x:=\{y\in \H \mid \proj_{\Gamma}y=x\}\] l'hyperplan orthogonal à $\Gamma$ en $x$. Nous définissons \textit{le tube $\Tube_{x,x'}^{\epsilon}$} comme l'enveloppe convexe de l'union des deux convexes fermés $\bar{B}(x,\epsilon)\cap \Gamma^{\perp}_x$ et $ \bar{B}(x',\epsilon)\cap \Gamma^{\perp}_{x'}$ (voir Figure \ref{figure_tube}, à noter que les tubes sont pleins). Le \textit{rayon} de ce tube en $z\in[x,x']$, noté $\rayon_{x,x'}^{\epsilon}(z)$, est \[\rayon_{x,x'}^{\epsilon}(z):=\sup\{\longueur{z}{u}\mid u\in \Gamma^{\perp}_z\cap \Tube_{x,x'}^{\epsilon}\}.\]
 Soient $x\leq y<y'\leq x'$ quatre points d'une même géodésique $\Gamma$. Nous dirons que le tube $\Tube_{x,x'}^{\epsilon_1}$ \textit{traverse} $\Tube_{y,y'}^{\epsilon_2}$ si nous avons les inégalités : $\rayon_{x,x'}^{\epsilon_1}(y)\leq \epsilon_2$ et $\rayon_{x,x'}^{\epsilon_1}(y')\leq \epsilon_2$ (voir Figure \ref{figure_tube_passant}). 	
		 		
		\begin{rmq}\label{rmq_convexite}
		Cette définition implique que pour tout $z\in[y,y']$, $\rayon_{x,x'}^{\epsilon_1}(z)\leq \rayon_{y,y'}^{\epsilon_2}(z)$. En effet, soit $u\in \Gamma^{\perp}_{z}\cap\Tube_{x,x'}^{\epsilon_1}$, alors $u$ appartient à un segment géodésique ayant ses extrémités dans $\Gamma^{\perp}_{y}\cap \Tube_{x,x'}^{\epsilon_1}\subset\Tube_{y,y'}^{\epsilon_2} $ et $\Gamma^{\perp}_{y'}\cap \Tube_{x,x'}^{\epsilon_1}\subset \Tube_{y,y'}^{\epsilon_2}$. Par convexité du tube $\Tube_{y,y'}^{\epsilon_2}$, $u$ appartient à $\Gamma^{\perp}_{z}\cap \Tube_{y,y'}^{\epsilon_2}$.
		\end{rmq}

	L'espace $\H$ étant $\cat(-1)$ et donc à fortiori $\cat(0)$, la distance est convexe. Cela revient à dire (voir \cite[p.176]{BH}) que pour tous segments géodésiques $[a,b]$ et $[c,d]$ et pour tous points $z\in[a,b]$ et $z'\in[c,d]$
		ayant les mêmes coordonnées barycentriques, c'est-à-dire $\longueur{a}{z}=t\longueur{a}{b}$ et $\longueur{c}{z'}=t\longueur{c}{d}$ avec $t\in[0,1]$, nous avons $\longueur{z}{z'}\leq (1-t)\longueur{a}{c}+t\longueur{b}{d}$. Nous utilisons la convexité de la distance de la façon suivante.
		
	\begin{lemme}\label{lemme_cvxd}
		Soit $G$ un groupe agissant sur $\H$. Considérons quatre points $x,x',z,z'$ sur une géodésique $\Gamma$ vérifiant $x\leq z<z'\leq x'$. Pour tout $\epsilon\geq 0$ nous avons : 
	\begin{enumerate}
		\item\label{lemme_cvxd_fix} $\Fix_{\epsilon}\{x,x'\}\subset \Fix_{\epsilon}\{z,z'\}$. 
		\item\label{lemme_cvxd_tube} le tube $\Tube_{x,x'}^{\epsilon}$ traverse le tube $\Tube_{z,z'}^{\epsilon}$.
	\end{enumerate}
	\end{lemme}	

	\begin{proof} Posons $t:=\frac{\dist(x,z)}{\dist(x,x')}$.
	Soit $g\in\Fix_{\epsilon}\{x,x'\}$. Comme $g$ agit par isométrie, $z\in[x,x']$ et $g\cdot z\in[g\cdot x,g\cdot x']$ ont les mêmes coordonnées barycentriques. Nous avons alors par convexité de la distance : \[\dist(z,g\cdot z)\leq t\dist(x',g\cdot x')+(1-t)\dist(x,g \cdot x)\leq\epsilon.\]
	En faisant de même pour $z'$, le point (\ref{lemme_cvxd_fix}) est démontré.

Par le même argument et en se plaçant dans un plan contenant le segment $[x,x']$, nous avons : \[\rayon_{x,x'}^{\epsilon}(z)\leq t\rayon_{x,x'}^{\epsilon}(x')+ (1-t)\rayon_{x,x'}^{\epsilon}(x) \leq\epsilon.\]
En faisant de même pour $z'$, le point (\ref{lemme_cvxd_tube}) est prouvé.
\end{proof}

		\subsection{Affaiblissement des hypothèses de la propriété WPD}\label{subsection_affaiblissement}
	Nous prouvons dans cette section le résultat suivant, qui nous permettra de vérifier plus facilement qu'une isométrie satisfait la propriété WPD.
		\begin{prop}\label{pro_equiv_WPD}
			Soient $G$ un groupe agissant sur $\H$ et $h$ un élément hyperbolique de $G$. Les deux propriétés suivantes sont équivalentes : \begin{enumerate}
				\item\label{propWPD_1} Il existe $y\in\Axe(h)$, $\eta>0$ et $n,k\in\N$ tels que $\Fix_{\eta}\{h^{-k}(y),h^{n}(y)\}$ soit fini.
				\item\label{propWPD_2} Il existe $w\in\Axe(h)$ tel que pour tout $\epsilon\geq 0$, il existe $N\in\N$ tel que $\Fix_{\epsilon}\{w,h^{N}(w)\}$ soit fini.
			\end{enumerate}
		\end{prop}	
		Le point (\ref{propWPD_2}) correspond à la définition donnée dans le lemme \ref{lemme_WPD_x}.(\ref{lemme_point_1}), et le point (\ref{propWPD_1}) est la version avec les hypothèses affaiblies qui nous sera utile dans la section \ref{section_cremona}. La preuve découle immédiatement du lemme \ref{lemme_cardfini}.
	\begin{lemme}\label{lemme_vois_tub}
	Soient $G$ un groupe agissant sur $\H$ et $\Gamma$ une géodésique de $\H$. Pour toutes constantes fixées $\epsilon\geq 0$, $\eta>0$ et pour tous points $y,y',z,z'\in\Gamma$ avec $y\leq z <z'\leq y'$, le tube $\Tube_{y,y'}^{\epsilon}$ traverse le tube $\Tube_{z,z'}^{\eta}$ dès que les distances $d(y,z)$ et $d(y',z')$ sont suffisamment grandes.
\end{lemme}

	\begin{figure}[h]
	 						\centering
	 						\def\svgwidth{220pt}
	 						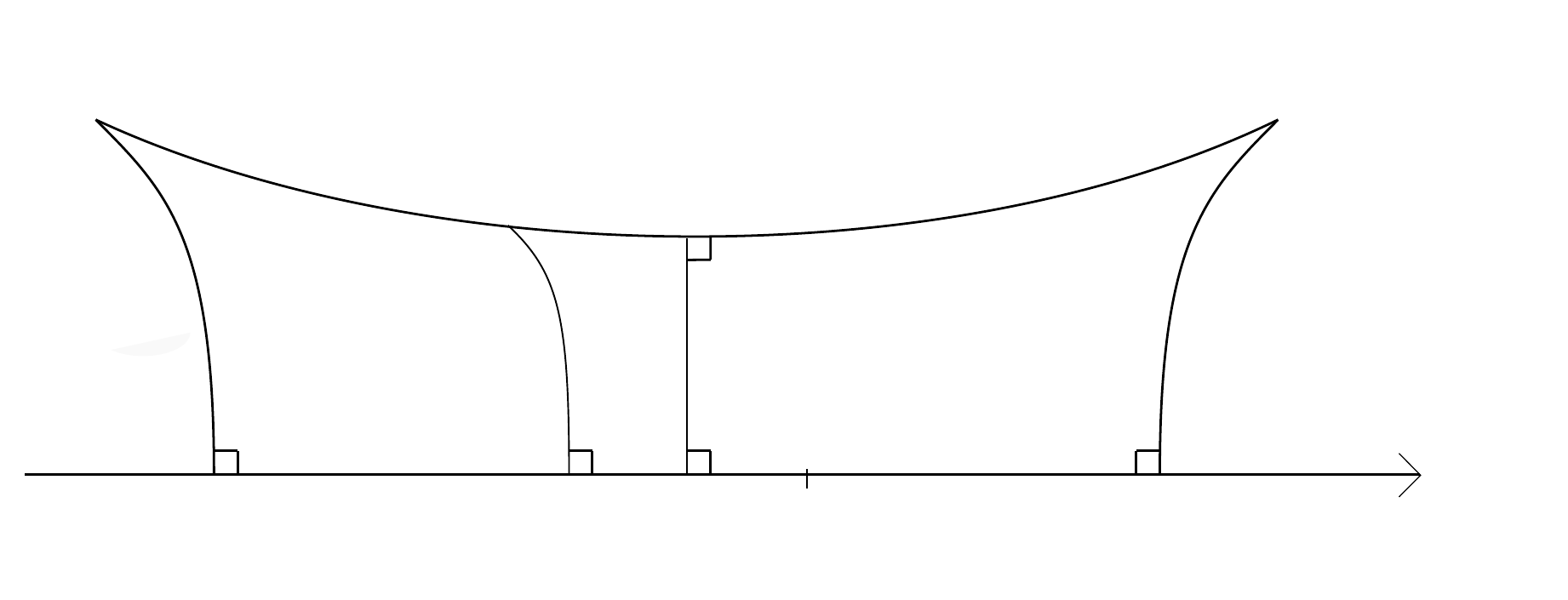
	 						\caption {Trace du tube $\Tube_{x,x'}^{\epsilon}$ dans le plan $\plan$.\label{figure_tube_plan}}
	 					\end{figure}

\begin{proof}
	Si $\eta\geq\epsilon$, il suffit de prendre $y\leq z$ et $y'\geq z'$ (voir lemme \ref{lemme_cvxd}.(\ref{lemme_cvxd_tube})).

	Si maintenant $\eta<\epsilon$, nous voulons trouver deux points $x$ et $x'$ sur $\Gamma$ avec $x\leq z$ et $x'\geq z'$ de sorte que $\rayon_{x,x'}^{\epsilon}(z)= \eta$ et $\rayon_{x,x'}^{\epsilon}(z')= \eta$. Ainsi en utilisant le lemme \ref{lemme_cvxd}.(\ref{lemme_cvxd_tube}), tous les points $y,y'$ vérifiant $d(y,z)\geq d(x,z)$ et $d(y',z')\geq d(x',z')$ satisferont le lemme. Soit $w$ le milieu de $[z,z']$. Considérons deux points $x,x'\in\Gamma$ avec $x\leq z$ et $x'\geq z'$ et tels que $w$ soit le milieu du segment $[x,x']$. Par symétrie, nous pouvons nous placer dans un plan $\plan$ contenant la géodésique $\Gamma$. Dans ce plan, la trace du tube $\Tube_{x,x'}^{\epsilon}$ est un quadrilatère plein. Notons $x_1$ et $x_1'$ deux sommets de ce quadrilatère situés dans un même demi-plan délimité par $\Gamma$ ({voir Figure \ref{figure_tube_plan}}). Nous avons alors $\longueur{x_1}{x}=\longueur{x_1'}{x'}=\epsilon$.
	 Soit $w_1$ le milieu du segment géodésique $[x_1,x_1']$. Par symétrie, la géodésique passant par les points $w$ et $w_1$ est la géodésique orthogonale aux géodésiques $[x_1,x_1']$ et $\Gamma$. Appelons $z_1$ le point sur le segment géodésique $[x_1,w_1]$ vérifiant $\proj_{\Gamma}z_1=z$.
	Le lemme \ref{lemme_relation-quadrilatere} nous donne : \[\left\{\begin{aligned}
	&\tanh\longueur{w}{w_1}\cosh\longueur{w}{z}=\tanh\longueur{z_1}{z}\\
	&\tanh\longueur{w}{w_1}\cosh\longueur{w}{x}=\tanh\longueur{x_1}{x}=\tanh\epsilon.
	\end{aligned}\right.\]
	En combinant les égalités, nous obtenons : \[\frac{\tanh\epsilon\cosh\longueur{w}{z}}{\cosh\longueur{w}{x}}=\tanh\longueur{z}{z_1}.\]
	Prenons $x$ tel que \[\longueur{w}{x}=\argcosh\left(\frac{\tanh\epsilon\cosh\longueur{w}{z}}{\tanh \eta}\right).\]
	Alors $\longueur{z}{z_1}=\eta$, et le point $w$ étant le milieu des segments $[z,z']$ et $[x,x']$, nous obtenons par symétrie que $\longueur{z'}{z'_1}=\eta$, comme attendu.
\end{proof}

	\begin{lemme}\label{lemme_cardfini}
	Soient $G$ un groupe agissant sur $\H$ et $h$ une isométrie hyperbolique de $G$. Pour toutes constantes fixées $\epsilon\geq 0$, $\eta >0$, pour tous points $z,z',w\in\Axe(h)$ avec $z<z'$, il existe $N,M\in\N$ tels que si l'ensemble $\Fix_{\eta}\{z,z'\}$ est fini alors $\Fix_\epsilon\{h^{-N}(w),h^M(w)\}$ l'est aussi.
	\end{lemme}
	
		\begin{figure}[h]
			\centering
			\def\svgwidth{300pt}
			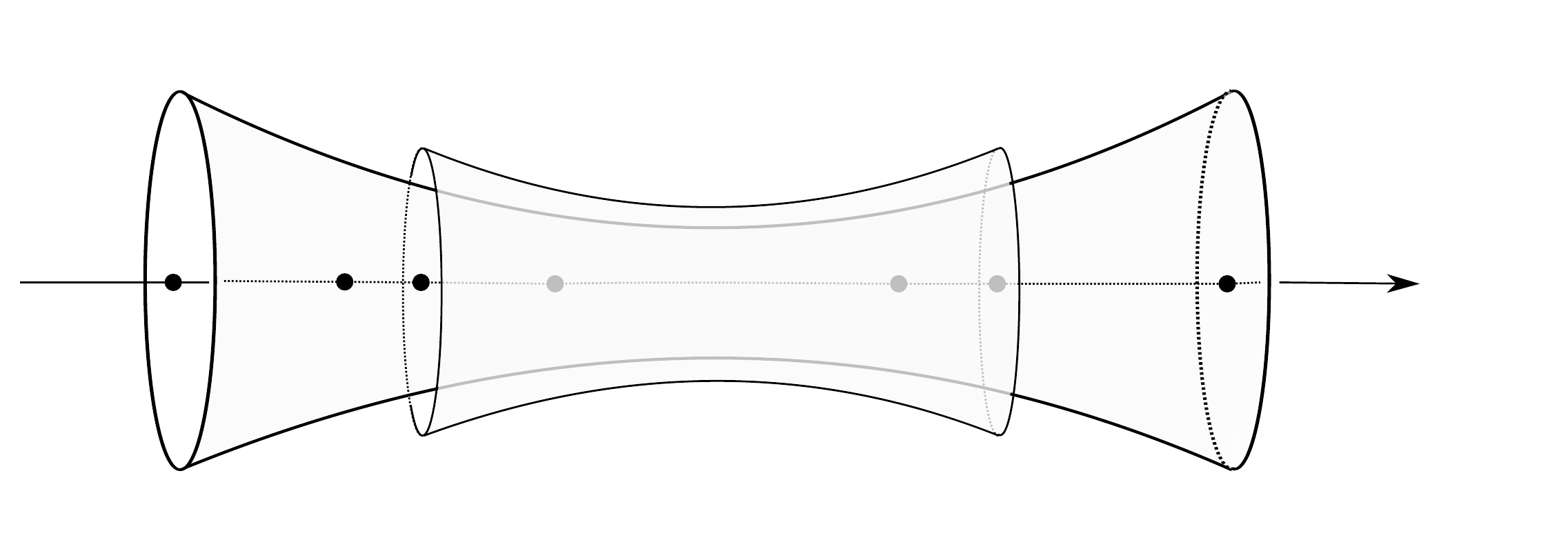
			\caption {Tube $\Tube_{h^{-N}(w)+\epsilon,h^M(w)-\epsilon}^{\epsilon}$ traversant le tube $\Tube_{z-\epsilon,z'+\epsilon}^{\eta/3}$.\label{figure_fix}}
		\end{figure}
	
	\begin{proof}
	Soient $\epsilon\geq0$, $\eta>0$, $z$,$z'$ et $w$ trois points de $\Axe(h)$. D'après le lemme \ref{lemme_vois_tub}, il existe deux entiers $N$ et $M$ suffisamment grands de sorte que le tube $\Tube_{h^{-N}(w)+\epsilon,h^M(w)-\epsilon}^{\epsilon}$ traverse $\Tube_{z-\epsilon,z'+\epsilon}^{\eta/3}$(voir Figure \ref{figure_fix}).

	Supposons que l'ensemble $\Fix_\epsilon\{h^{-N}(w),h^M(w)\}$ soit infini. Considérons une suite $(f_n)_{n\in \N}$ d'éléments deux à deux distincts de cet ensemble, et notons $p_n$ et $p_n'$ les projections respectives de $f_n(h^{-N}(w))$ et $f_n(h^M(w))$ sur l'axe de $h$. Ainsi, pour tout $n$, nous avons : \[[f_n(h^{-N}(w)), f_n(h^M(w))]\subset \Tube_{p_n,p_n'}^{\epsilon}.\] En appliquant le lemme \ref{lemme_cvxd}.(\ref{lemme_cvxd_tube}) nous obtenons que le tube $\Tube_{p_n,p_n'}^{\epsilon}$ traverse le tube $\Tube_{h^{-N}(w)+\epsilon,h^M(w)-\epsilon}^{\epsilon}$ qui, par construction, traverse le tube $\Tube_{z-\epsilon,z'+\epsilon}^{\eta/3}$.
	De plus, par le lemme \ref{lemme_cvxd}.(\ref{lemme_cvxd_fix}), nous obtenons l'inclusion : \[\Fix_\epsilon\{h^{-N}(w),h^M(w)\}\subset \Fix_\epsilon\{z,z'\}.\] Par conséquent, pour tout $n$, $f_n(z)$ appartient à la boule fermée $\bar{B}(z,\epsilon)$. Notons $z_n$ la projection de $f_n(z)$ sur $\Axe(h)$. En particulier, la suite $(z_n)_{n\in\N}$ appartient au compact $[z-\epsilon,z'+\epsilon]$. Ceci signifie que d'une part, quitte à prendre une sous-suite, nous pouvons supposer que la suite $(z_n)_{n\in\N}$ converge et que pour tout $k\geq 0$, $\longueur{z_0}{z_k}\leq\frac{\eta}{3}$, et d'autre part que $\dist(z_n,f_n(z))\leq\frac{\eta}{3}$.
	Par conséquent, nous avons pour tout $n\geq 0$ : \[\dist(f_0(z),f_n(z))\leq \dist(f_0(z),z_0)+ \dist(z_0,z_n) +\dist(z_n,f_n(z))\leq\eta.\]
	Quitte à prendre une sous-suite de $(f_n)_{n\in\N}$ et en considérant maintenant le point $z'$, nous obtenons, par le même argument, l'inégalité précédente en remplaçant $z$ par $z'$.	
	Ainsi, la suite $(f_0^{-1}f_n)_{n\geq 0}$ est contenue dans $\Fix_{\eta}\{z,z'\}$ qui est donc infini comme attendu.
	\end{proof}

	\section{Application au groupe de Cremona}\label{section_cremona}
	
	Soient $\K$ un corps et $\bar{\K}$ sa clôture algébrique. Dans cette section, nous rappelons la construction de l'espace de Picard-Manin sur $\bar{\K}$, et celle de l'espace hyperbolique de dimension infinie associé. Ainsi, nous pourrons montrer que le groupe $\Bir(\P_{\K})$ agit discrètement le long de l'axe de certains de ses éléments en utilisant la version affaiblie établie dans la section \ref{subsection_affaiblissement}. 
	\subsection{Action du groupe de Cremona sur l'espace de Picard-Manin}\label{subsection_picard_manin}
	Plaçons nous sur $\bar{\K}$. Rappelons rapidement la construction de l'espace de Picard-Manin et comment le groupe $\Bir(\P_{\bar{\K}})$ agit dessus. Plus de précisions se trouvent dans \cite{CL} et \cite{C}. 
	
	Soit $S$ une surface projective lisse. Le groupe de Néron-Severi associé à $S$, noté $\NS(S)$, est le groupe des diviseurs à coefficients réels sur $S$ à équivalence numérique près. Il est muni d'une forme bilinéaire symétrique, la forme d'intersection. Pour tout diviseur $D$ sur $S$ nous notons $\{D\}_S$ sa classe de Néron-Severi ou $\{D\}$ s'il n'y a pas d’ambiguïté sur la surface. Si $\pi : S'\longrightarrow S$ est un morphisme birationnel entre deux surfaces, alors le tiré en arrière \[\pi^* : \NS(S)\hookrightarrow \NS(S')\] qui à la classe d'un diviseur associe la classe de sa transformée totale, est un morphisme injectif qui préserve la forme d'intersection. De plus, $\NS(S')$ est isomorphe à
	\begin{equation}\label{eq_NS}
	\pi^*(\NS(S))\oplus(\underset{p\in\Bs(\pi^{-1})}{\oplus}\R\{E_p^*\}),
	\end{equation}
	 où $\Bs(\pi^{-1})$ est l'ensemble des points base de $\pi^{-1}$ (infiniment proches ou pas) et $E_p^*$ est la transformée totale, vue dans $S'$, du diviseur exceptionnel $E_p$ obtenu en éclatant le point $p$. Cette somme est orthogonale relativement à la forme d'intersection.

	Considérons la limite inductive des groupes de Néron-Severi des surfaces $S'$ obtenues comme éclatements de $S$ (nous dirons que de telles surfaces $S'$ dominent $S$) : \[\PM_C(S)=\lim\limits_{\underset{S'\rightarrow S}{\longrightarrow}}\NS(S'),\] 
	où l'indice $C$ fait référence aux b-diviseurs de Cartier (pour plus de précisions voir \cite{Fa}).
	Remarquons que pour toute surface $S'$ dominant $S$, le groupe $\NS(S')$ est plongé dans $\PM_C(S)$.
	En fait, si nous considérons un diviseur $D$ sur $S$, à chaque surface $S'$ dominant $S$, nous pouvons lui faire correspondre une classe de Néron-Severi $\{D \}_{S'}$ dans $\NS(S')$. Ces éléments sont tous identifiés dans $\PM_C(S)$ et correspondent à une classe $d$ de $\PM_C(S)$ notée en lettre minuscule.
	\begin{ex}\label{ex_clasNS}
		Considérons un point $q$ appartenant au diviseur exceptionnel $E_p$, issu de l'éclatement d'une surface $S$ au point $p$. Notons $S_p$ la surface obtenue en éclatant le point $p$ et $S_{p,q}$ celle en éclatant successivement les points $p$ et $q$. La classe $e_p$ correspond à $\{E_p\}_{S_p}$ dans $\NS(S_p)$ et à $\{\tilde{E}_p+E_q\}_{S_{p,q}}$ dans $\NS(S_{p,q})$ où $\tilde{E}_p$ est la transformée stricte de $E_p$ dans $S_{p,q}$.
	\end{ex}
	 Définissons la forme d'intersection sur $\PM_C(S)$. Pour cela, considérons $c$ et $d$ deux éléments de $\PM_C(S)$. Il existe une surface $S_1$ dominant $S$ telles que les classes $c$ et $d$ correspondent respectivement à $\{C\}_{S_1}$ et $\{D\}_{S_1}$ dans $\NS(S_1)$. La forme d'intersection est donnée par : $c \cdot d=\{C\}_{S_1}\cdot \{D\}_{S_1}$. Montrons qu'elle ne dépend pas du choix de la surface $S_1$. En effet, considérons une autre surface $S_2$ où les éléments $c$ et $d$ se réalisent comme des éléments de $\NS(S_2)$. En résolvant l'application birationnelle allant de $S_1$ vers $S_2$ nous obtenons une surface $S_3$ dominant $S_1$ et $S_2$ c'est-à-dire qu'il existe deux morphismes birationnels $\pi :S_3\rightarrow S_1$ et $\sigma : S_3\rightarrow S_2$ faisant commuter le diagramme : 	\begin{center}
	 		\begin{tikzcd}[column sep=small]
	 			& S_3\arrow{dl}[above left]{\pi} \arrow{dr}{\sigma} & \\
	 			S_1 \arrow[dashrightarrow]{rr} \arrow{dr} & & S_2 \arrow{dl}.\\
	 			& S &
	 		\end{tikzcd}		
	 	\end{center}
	La forme d'intersection étant stable par tiré en arrière nous avons : \[\{C\}_{S_1}\cdot \{D\}_{S_1}=\{C\}_{S_3}\cdot \{D\}_{S_3}=\{C\}_{S_2}\cdot \{D\}_{S_2}.\]

	Par la suite, nous nous intéressons à l'espace de Hilbert défini par  \[\PM(S)=\{\{D_0\}_{S}+\sum\limits_{p}\lambda_pe_p\mid \lambda_p\in\R,\ \sum\lambda_p^2<\infty \text{ et }\{D_0\}_{S} \in\NS(S)\},\]
	que nous appelons l'espace de Picard-Manin (voir \cite{CL} et \cite{C}). C'est le complété $L^2$ de $\PM_C(S)$. Les classes $e_p$ (nous gardons les notations introduites dans l'exemple \ref{ex_clasNS}) où $p$ est un point de $S$ ou d'une surface dominant $S$, sont d'auto-intersection $-1$, orthogonales deux à deux et orthogonales à $\NS(S)$. La forme d'intersection est donc de signature $(1,\infty)$ et préserve la décomposition orthogonale :
	\begin{equation}\label{eq_PM}
	\PM(S)=\NS(S)\oplus (\underset{p}{\oplus}\R e_p).
	\end{equation}

	\begin{ex}
	Soient $p$ et $q$ les points définis dans l'exemple \ref{ex_clasNS}. Considérons les classes de Picard-Manin $e_p$ et $e_q$. Plaçons nous sur la surface $S_{p,q}$ puisqu'elle domine la surface $S_p$ contenant le diviseur $E_p$ et qu'elle contient le diviseur $E_q$. Dans $\NS(S_{p,q})$ les classes $e_p$ et $e_q$ correspondent respectivement à $\{\tilde{E}_p+E_q\}$ et à $\{E_q\}$. Nous avons donc : \begin{align*}
	& e_p\cdot e_p = \tilde{E}_p^2+E_q^2+2\tilde{E}_p\cdot E_q=-2-1+2=-1\\
	& e_p\cdot e_q = \tilde{E}_p\cdot E_q +E_q^2=1-1=0. 
	\end{align*} 
	\end{ex}

 	Si $\pi:S'\rightarrow S$ est un morphisme birationnel, nous avons mentionné que l'application induite par tiré en arrière sur les groupes de Néron-Severi est une injection. Précisément : \[\begin{array}{lccc}
 	\pi^*: & \NS(S) & \hookrightarrow & \NS(S')\\
 	& \{D_0\}_{S} & \mapsto & \{\tilde{D_0}\}_{S'}+\sum\limits_{q\in\Bs(\pi^{-1})}m_q(D_0)e_q,
 	\end{array}\]
 	où $\tilde{D}_0$ est la transformée stricte dans $S'$ de $D_0$ et $m_q(D_0)$ la multiplicité de $D_0$ au point $q$. Regardons ce que nous obtenons au niveau des espaces de Picard-Manin. D'après (\ref{eq_PM}) puis (\ref{eq_NS}), nous avons la décomposition \[\PM(S')=\left(\NS(S)\oplus (\underset{q\in \Bs(\pi^{-1})}{\oplus}\R e_q)\right)\oplus (\underset{r\notin \Bs(\pi^{-1})}{\oplus}\R e_r).\] Nous définissons un isomorphisme $\isome{\pi}{}{}$ de $\PM(S')$ vers $\PM(S)$ qui consiste à considérer comme exceptionnelles au-dessus de $S$ les classes $e_q$ qui étaient dans l'espace de Néron-Severi de $S'$. 
 	Ainsi l'injection $\pi^*$ au niveau des espaces de Néron-Severi devient un isomorphisme $\isome{\pi}{}{-1}$ au niveau des espaces de Picard-Manin.
 	 Plus précisément, nous avons \[\begin{array}{lccc}
	\isome{\pi}{}{-1} : & \PM(S)=\NS(S)\oplus (\underset{p}{\oplus}\R e_p) & \longrightarrow & \PM(S')=\NS(S')\oplus (\underset{r\notin\Bs(\pi^{-1})}{\oplus}\R e_r)\\
	& \{D_0\}_{S}+\sum\limits_{p}\lambda_pe_p  & \mapsto &  \left(\{\tilde{D_0}\}_{S'}+\sum\limits_{q\in\Bs(\pi^{-1})}(m_q(D_0)+\lambda_q)e_q\right)+\sum\limits_{r\notin\Bs(\pi^{-1})}\lambda_re_r.
	\end{array}\] 
	
	À présent, considérons l'espace \[\H(S)=\{c\in \PM(S)\mid c\cdot c=1 \text{ et } c\cdot d_0>0\},\] où $d_0\in\NS(S)$ est une classe ample. Muni de la distance définie par $\dist(c,c')=\argcosh(c\cdot c')$ pour tous $c,c'\in \H(S)$, c'est un espace hyperbolique de dimension infinie comme ceux introduits dans la section \ref{section_hyperbolique}. Nous nous intéressons dans cet article à $\H(\P_{\bar{\K}})$ que nous notons $\H_{\bar{\K}}$. Tout élément de $\H_{\bar{\K}}$ est de la forme $\lambda_{\ell}\ell +\sum_{p}\lambda_pe_p$ où $\ell$ correspond à la classe d'une droite dans $\P_{\bar{\K}}$, $\lambda_{\ell}>0$ et $\lambda_{\ell}^2-\sum_{p}\lambda_p^2=1$.
	Considérons une résolution de $f\in\Bir(\P_{\bar{\K}})$ :
	\begin{center}
		\begin{tikzcd}[column sep=small]
			& S\arrow{dl}[above left]{\pi} \arrow{dr}{\sigma} & \\
			\P_{\bar{\K}} \arrow[dashrightarrow]{rr}{f} & & \P_{\bar{\K}}.
		\end{tikzcd}		
	\end{center}
	Le groupe $\Bir(\P_{\bar{\K}})$ agit sur $\H_{\bar{\K}}$ via l'application $(f,c)\mapsto \isome{f}{}{}(c)$ où $\isome{f}{}{}$ est définie par \[\isome{f}{}{}=\isome{\sigma}{}{}\circ (\isome{\pi}{}{})^{-1}.\] Remarquons que $(\isome{f}{}{})^{-1}=\isome{(f^{-1})}{}{}$.

		\begin{rmq}\label{rmq_isom_local}
			Si $f$ est un isomorphisme d'un voisinage $U$ de $x$ sur un voisinage $V$ de $f(x)$ alors nous avons
			$\isome{f}{}{}(e_x)=e_{f(x)}$.
		\end{rmq}

		Explicitons sur un exemple l'action de $f$ sur $\ell$, la classe d'une droite de $\P_{\bar{k}}$.

	\begin{ex}
		Soit $f : \P_s\dashrightarrow \P_b$ une application quadratique et notons $q_1$, $q_2$ et $q_3$ les points base (infiniment proches ou non) de $f^{-1}$. Nous mettons un indice $s$ comme source et un indice $b$ comme but pour plus de clarté.
		La transformée par $f$ d'une droite générale $L$ de $\P_s$ (ne passant pas par les points base de $f$) est une conique $C$ de $\P_b$ qui passe par les points base de $f^{-1}$.
		Considérons une résolution $\pi : S \rightarrow \P_s$ et $\sigma : S \rightarrow \P_b$ de $f$. Nous voulons exprimer $\{\pi^*(L)\}_s$ dans la base de $\NS(S)$ issue de $\sigma^*$, c'est-à-dire dans la base : $\sigma^*(\{L\}_b)$, $\{E^*_{q_1}\}_b$, $\{E^*_{q_2}\}_b$, $\{E^*_{q_3}\}_b$. La droite $L$ ne passant pas par les points base de $f$ nous avons $\pi^*(L)=\tilde{L}$ qui correspond à la transformée stricte de $C$ par $\sigma$ : $\tilde{C}=\sigma^*(C)-E^*_{q_1}-E^*_{q_2}-E^*_{q_3}$. Or, $\sigma^*(\{C\}_b)=2\sigma^*(\{L\}_b)$. Ainsi, nous avons écrit $\{\pi^*(L)\}_s$ dans la base voulue, $\{\pi^*(L)\}_s=\{2\sigma^*(L)-E^*_{q_1}-E^*_{q_2}-E^*_{q_3}\}_b$. Si nous réécrivons cette égalité dans les espaces de Picard-Manin nous avons : \[\isome{f}{}{}(\ell)=2\ell-e_{q_1}-e_{q_2}-e_{q_3}.\]

	\end{ex}
	De manière générale, pour tout $f\in\Bir(\P_{\bar{\K}})$, \[\isome{f}{}{}(\ell)\cdot\ell=\deg(f).\]

	\subsection{Action des applications $h_n$ sur $\H_{\bar{k}}$ }
	Par la suite, nous allons nous intéresser aux applications $h_n:(x,y)\mapsto(y,y^n-x)$ où $n\geq 2$.
	 Elles appartiennent à $\Bir(\P_{\K})$ qui est inclus dans $\Bir(\P_{\bar{k}})$. Nous nous plaçons donc sur $\bar{\K}$. Fixons l'entier $n\geq 2$. 
	 
\begin{figure}
			\centering
			\def\svgwidth{200pt}
			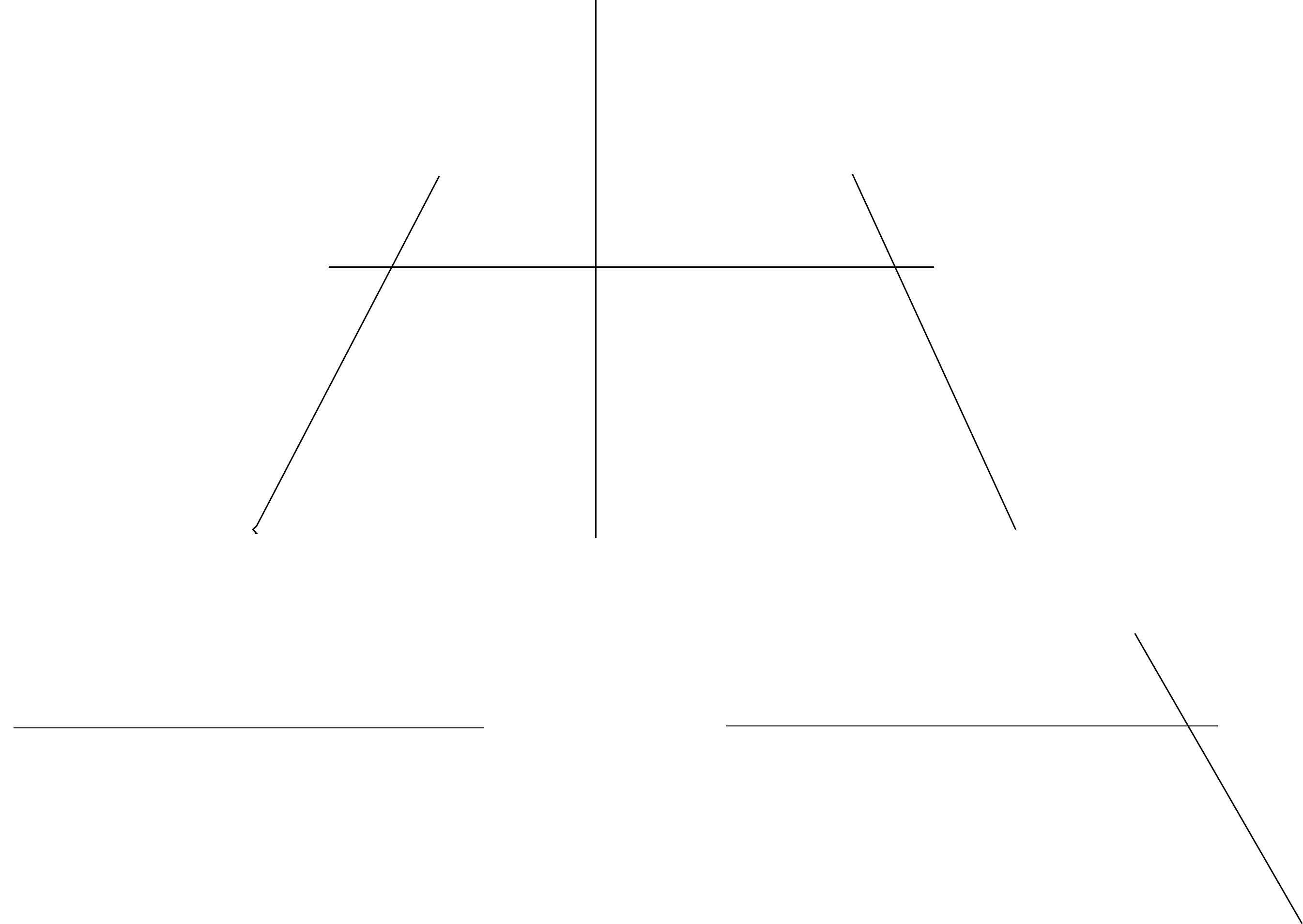
			\caption{Arbre des transformées strictes des diviseurs exceptionnels obtenus en résolvant les points d'indétermination de $h_n$.\label{fig_arbres}}
		\end{figure}

	 Considérons les applications birationnelles $a: (x,y)\mapsto (y,x)$ et $j_n : (x,y)\mapsto (y^n-x,y)$. Nous avons : \[h_n=a\circ j_n \text{ et } h_n^{-1}=a^{-1}\circ h_n\circ a.\] L'application $j_n$ étant une application de Jonquières elle a un point base $p_0$ de multiplicité $n-1$ et $2n-2$ points base de multiplicité $1$ (voir \cite[Definition 2.28]{KS}). En homogénéisant, le point $p_0$ a pour coordonnées $[1:0:0]$. C'est le seul point base de l'application $j_n$ qui appartient à $\P_{\bar{k}}$. Les autres points base de $j_n$ forment une tour de points infiniment proches au-dessus de $p_0$ (voir \cite[Lemme 9]{La}). Pour $1\leq k\leq 2n-2$, notons $p_k$ le point base de $j_n$ appartenant au diviseur exceptionnel obtenu en éclatant le point base $p_{k-1}$. Ainsi le point $p_k$ est infiniment proche du point $p_{k-1}$. Plus précisément, la figure \ref{fig_arbres} représente les transformées strictes des diviseurs exceptionnels obtenus en résolvant $h_n$ munis de leur auto-intersection respective. La configuration obtenue, illustrée sur la figure, découle du fait que les $n-1$ points base $p_1,\dots,p_{n-1}$ se trouvent sur les transformées strictes respectives de $E_{p_0}$. 
	À l'exception des points $p_0$ et $p_1$, remarquons que les points $p_k$ dépendent de $n$. Cependant, pour ne pas alourdir davantage l'écriture nous ne mettons pas d'indice $n$.
	L'application $a$ étant un automorphisme de $\P_{\bar{\K}}$, $h_n$ possède les mêmes points base que $j_n$. En conjuguant par $a$, nous obtenons par symétrie les points base de $h_n^{-1}$ que nous notons $\{q_k\}_{0\leq k \leq 2n-2}$. Remarquons que $q_0=[0:1:0]$.

	L'isométrie $\isome{h}{n}{}$ est hyperbolique et de longueur de translation $\log(n)$ (voir \cite[Theorem 4.4, Remark 4.5]{CL}). Nous allons décrire son axe.
	Pour alléger l'écriture nous notons $e_n^+$ (respectivement $e_n^-$) les sommes avec multiplicités des classes des diviseurs exceptionnels issus de la résolution de $h_n$ (respectivement de $h_n^{-1}$) : \[\left\{ \begin{aligned}
	e_n^+&=(n-1)e_{p_0} + e_{p_1} +\dots + e_{p_{2n-2}}\\
	e_n^-&=(n-1)e_{q_0} + e_{q_1} +\dots + e_{q_{2n-2}}
	\end{aligned}\right. .\]
	L'action sur $\ell$ de $\isome{h}{n}{}$ et de ses itérés est donnée par : \[\isome{h}{n}{}(\ell)=n\ell-e_n^-\text{, } \isome{h}{n}{2}(\ell)=n^2\ell-ne_n^--\isome{h}{n}{}(e_n^-),\text{  etc.}\]
	La suite $(\frac{1}{n^k}\isome{h}{n}{k}(\ell))_{k\in\N}$ converge dans l'espace de Picard-Manin vers un élément $b_n^+$ d'auto-intersection $0$ qui s'identifie à un point du bord $\partial \H_{\bar{\K}}$. Cet élément correspond à un bout de $\Axe(\isome{h}{n}{})$.
	De même, la suite $(\frac{1}{n^k}\isome{h}{n}{-k}(\ell))_{k\in\N}$ converge vers un élément $b_n^-$. Ces deux classes se réécrivent :\[b_n^+=\ell-\sum\limits_{i=0}^{\infty}\frac{\isome{h}{n}{i}(e_n^-)}{n^{i+1}}\text{ et } b_n^-=\ell-\sum\limits_{i=0}^{\infty}\frac{\isome{h}{n}{-i}(e_n^+)}{n^{i+1}}.\]
	
	L'action de $h_n$ sur les classes $e_p$ est donnée par la remarque suivante car $h_n$ est un automorphisme polynomial.
	\begin{rmq}\label{rmq_heq}
	Soit $f$ un automorphisme du plan affine $\A^2$ de degré au moins $2$, que nous étendons en une application birationnelle de $\P$.
	L'ordre des éclatements dans la résolution minimale de $f^{-1}$ est uniquement déterminé. En effet à chaque étape de la résolution l'unique point base propre est l'image de la droite à l'infini.
	Ainsi il y a dans ce contexte une notion de \og dernier diviseur \fg \ produit lors de la résolution minimale de $f^{-1}$, et si $p$ est un point de la droite à l'infini qui n'est pas le point base propre de $f$, alors il existe un point $q$ sur ce dernier diviseur produit tel que la transformation birationnelle induite par $f$ donne un isomorphisme local entre des voisinages de $p$ et $q$.
	En particulier, par la remarque \ref{rmq_isom_local}, nous avons $\isome{f}{}{}(e_p) = e_q$. 
	\end{rmq}

			\begin{lemme}\label{lemme_tousdistincts}
			Considérons l'ensemble des éléments des $4n-2$ suites $(\isome{h}{n}{i}(e_{q_k}))_{i\in\N}$, $(\isome{h}{n}{-i}(e_{p_k}))_{i\in\N}$ avec $0\leq k\leq 2n-2$. Ces éléments sont tous deux à deux orthogonaux.
			\end{lemme}

				\begin{proof}
				Construisons les points $(\isome{h}{n}{}(e_{q_k}))$, pour $0\leq k\leq 2n-2$. Les autres points se construisent de la même façon.
				L'application $h_n$ est régulière en $q_0$ et $h_n(q_0)=q_0$ est un point base de $h_n^{-1}$. D'après la remarque \ref{rmq_heq}, il existe un point sur le dernier diviseur exceptionnel $E_{q_{2n-2}}$ issu de l'éclatement de tous les points $\{q_k\}_{0\leq k\leq 2n-2}$, que nous notons $q_{2n-1}$, tel que $\isome{h}{n}{}(e_{q_0})=e_{q_{2n-1}}$. Si nous notons $\pi_0$ l'éclatement du point $q_0$, alors l'application $h_n\circ \pi_0$ est régulière en $q_1$ et $h_n\circ \pi_0(q_1)=q_0$ qui est un point base de $(h_n\circ \pi_0)^{-1}$. Toujours d'après la remarque \ref{rmq_heq}, il existe un point sur le diviseur exceptionnel obtenu en éclatant le point $q_{2n-1}$, que nous notons $q_{2n}$, et vérifiant $\isome{h}{n}{}(e_{q_1})=e_{q_{2n}}$.
				En réitérant le procédé, nous obtenons qu'il existe un point, noté $q_{2n-1+k}$, sur le diviseur exceptionnel issu de l'éclatement du point $q_{2n-2+k}$ tel que $\isome{h}{n}{}(e_{q_k})=e_{q_{2n-1+k}}$.
				
				À présent, nous pouvons montrer, par récurrence sur $i$, que $\isome{h}{n}{i}(e_{q_k})=e_{q_{i(2n-1)+k}}$ pour tout $i\geq 1$ et tout $k$ tel que $0\leq k \leq 2n-2$.

				En faisant de même, nous obtenons également que $\isome{h}{n}{-i}(e_{p_k})=e_{p_{i(2n-1)+k}}$ pour tout $i\geq 1$ et tout $k$ tel que $0\leq k \leq 2n-2$. Ainsi tous les éléments de toutes les suites sont deux à deux orthogonaux.
				\end{proof}

	Considérons à présent pour chaque $n\geq 2$, le point $w_n$ qui est le projeté de $\ell$ sur $\Axe(\isome{h}{n}{})$. 
	L'axe de $\isome{h}{n}{}$ étant uniquement déterminé par $b_n^+$ et $b_n^-$, $w_n$ est une combinaison linéaire de ces deux classes; $w_n=\alpha b_n^++\beta b_n^-$. Nous avons $1=w_n^2=2\alpha \beta$ car $(b_n^+)^2=0=(b_n^-)^2$ et $b_n^+\cdot b_n^-=1$. De plus, $w_n\cdot \ell=\alpha +\beta$ doit être minimal car $w_n$ est le projeté de $\ell$. Nous obtenons finalement : \begin{equation}\label{eq_wn}
	w_n=\sqrt{2}\ell-\frac{1}{\sqrt{2}}r_n \text{ \  où \  } r_n=\sum\limits_{i=0}^{\infty}\frac{\isome{h}{n}{i}(e_n^-)+\isome{h}{n}{-i}(e_n^+)}{n^{i+1}}.
	\end{equation}

		\begin{rmq}\label{rmq_tousdistincts}
		Écrivons tous les termes de $r_n$ : \[r_n=(n-1)\dfrac{e_{q_0}}{n}+ \dfrac{e_{q_1}}{n} +\dots + \dfrac{e_{q_{2n-2}}}{n}+(n-1)\dfrac{e_{p_0}}{n}+ \dots + \dfrac{e_{p_{2n-2}}}{n} + (n-1)\dfrac{\isome{h}{n}{}(e_{q_0})}{n^2}+\cdots .\]
		Le lemme \ref{lemme_tousdistincts} implique que tous les termes de $r_n$ sont orthogonaux deux à deux. Ainsi, la classe de tout diviseur exceptionnel $e_i$ a un nombre d'intersection non nul avec au plus un seul des termes de $r_n$.
		
		\end{rmq}
		
		Remarquons également, même si nous n'utiliserons pas ce fait, que les points $p_i$ (respectivement $q_i$), construits dans le lemme \ref{lemme_tousdistincts}, sont les points base des itérés de $h_n$ (respectivement de $h_n^{-1}$).
		Le diagramme suivant donne une idée de preuve dans le cas de $h_n^2$.

			\begin{center}
				\begin{tikzcd}[column sep=small]
					& & \tilde{S} \arrow{dl}[above left]{\tilde{\pi}}  \arrow{dr}{\tilde{\sigma}} & &\\
					& S\arrow{dl}[above left]{\pi} \arrow{dr}{\sigma} & & S\arrow{dl}[above left]{\pi} \arrow{dr}{\sigma} &\\
					\P_{\bar{\K}} \arrow[dashrightarrow]{rr}{h_n} & & \P_{\bar{\K}} \arrow[dashrightarrow]{rr}{h_n} & & \P_{\bar{\K}}
				\end{tikzcd}		
			\end{center}

	\subsection{Preuve du résultat principal}\label{section_resprincipal}
	Soient ${\K}$ un corps et $n\geq  2$ un entier tels que la caractéristique de ${\K}$ ne divise pas $n$. Considérons les applications \[\begin{array}[t]{lrcl}
				h_n : & \A_{\K}^2 & \longrightarrow & \A_{\K}^2 \\
				& (x,y) & \longmapsto & (y,y^n-x) \end{array}\] appartenant au groupe $\Bir(\P_{\K})$. 
	Nous montrons dans ce paragraphe le théorème \ref{prop_principale}, c'est-à-dire que le groupe $\Bir(\P_{\K})$ agit discrètement le long de l'axe de $\isome{h}{n}{}$. Plus précisément, nous allons montrer que le groupe $\Bir(\P_{\bar{\K}})$ agit discrètement le long de l'axe de $\isome{h}{n}{}$, ce qui impliquera immédiatement la même propriété pour le groupe $\Bir(\P_{\K})$.
	Nous avons besoin du lemme suivant tiré de \cite[Proposition 5.7]{CL}. Remarquons que nous avons effectué deux modifications dans notre version de l'énoncé.
	D'une part nous voulons vérifier le lemme pour les transformations $h_n$, qui ne sont pas des applications générales.
	D'autre part le rapporteur, que nous remercions, nous a signalé une erreur de calcul d'intersection dans la preuve, qui en fait se trouvait déjà dans \cite{CL}.
	C'est ce qui explique que dans notre énoncé la constante epsilon dépende de $n$, alors qu'elle était uniforme dans l'article original.
	Ceci est suffisant pour notre but, signalons cependant qu'il serait possible de préserver l'uniformité de la constante au prix d'une preuve un peu plus longue. 
	
	Avant d'énoncer le lemme, rappelons que $w_n$ est le projeté sur l'axe de $\isome{h}{n}{}$ de la classe de $\ell$ et est donné par l'égalité (\ref{eq_wn}).
	
		\begin{lemme}\label{lemme_lineaire}
Il existe $\epsilon_n>0$ tel que tout élément de $\Fix_{\epsilon_n}\{w_n\}\subset\Bir(\P_{\bar{\K}})$ soit un automorphisme de $\P_{\bar{\K}}$.
		\end{lemme}

		\begin{proof}
		Choisissons $\epsilon_n\in \left]0, \argcosh\Big(\sqrt{2}+\frac{1}{n\sqrt{2}}\Big)-\argcosh\sqrt{2}\right]$, cette constante vérifie alors : \[\label{hypothese_star}\tag{$\ast$}\left\{\begin{aligned}
		&\argcosh(\sqrt{2})+\epsilon_n<\argcosh\Big(\sqrt{2}+\frac{1}{n\sqrt{2}}\Big)<\argcosh\frac{3}{\sqrt{2}}\\
		&2\argcosh(\sqrt{2})+\epsilon_n<\argcosh4
		\end{aligned}\right. .\]
		La deuxième inégalité vient du fait que :\[2\argcosh(\sqrt{2})+\epsilon_n<\argcosh(\sqrt{2})+\argcosh\Big(\sqrt{2}+\frac{1}{n\sqrt{2}}\Big)<\argcosh(\sqrt{2})+\argcosh\Big(\frac{5}{2\sqrt{2}}\Big)<\argcosh 4,\]
		avec $\argcosh(\sqrt{2})+\argcosh\Big(\frac{5}{2\sqrt{2}}\Big)\simeq 0,881+1,171=2,052$ et $\argcosh 4 \simeq 2,063$.

		Soit $f\in\Fix_{\epsilon_n}\{w_n\}$. Nous voulons montrer que $f$ est de degré $1$. D'après l'inégalité triangulaire et le fait que $\dist(\isome{f}{}{}(w_n),w_n)\leq \epsilon_n$, nous obtenons (voir Figure \ref{figure_projection}) : \begin{align*}
		\dist(\isome{f}{}{}(\ell),\ell)&\leq\dist(\isome{f}{}{}(\ell),\isome{f}{}{}(w_n))+\dist(\isome{f}{}{}(w_n),w_n)+\dist(w_n,\ell)\\
		&\leq 2\dist(w_n,\ell)+\epsilon_n.
		\end{align*}
		
		\begin{figure}
				\centering
				\def\svgwidth{140pt}
				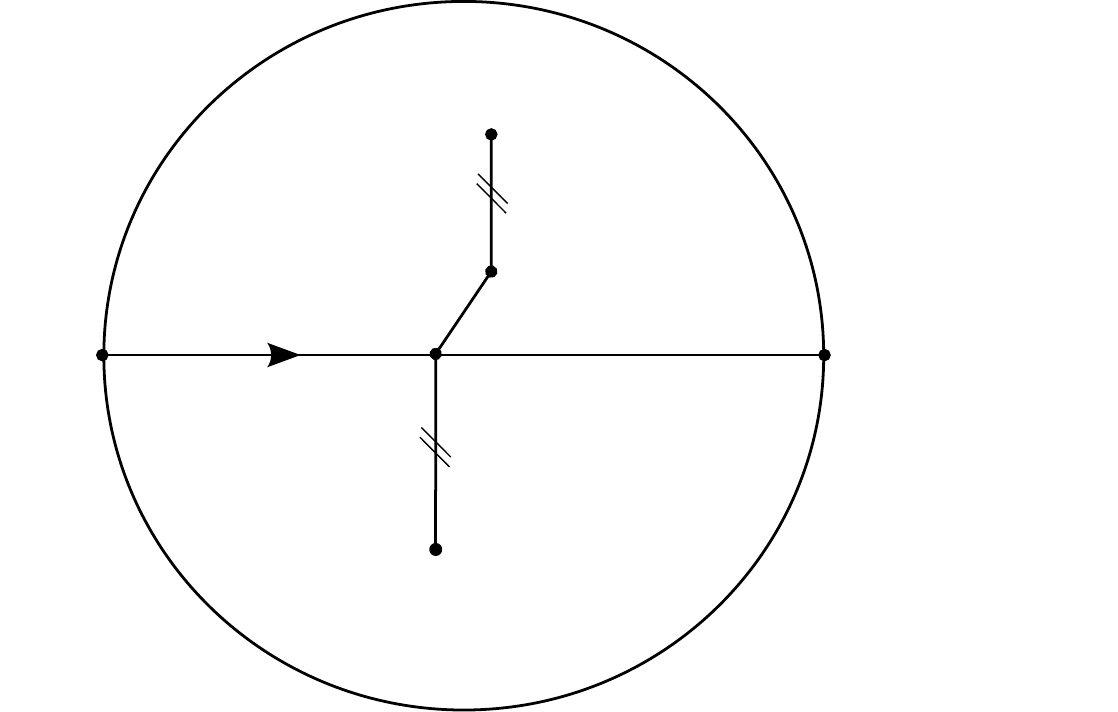
				\caption{Action de $f$ sur les classes $\ell$ et $w_n$. \label{figure_projection}}
			\end{figure}
		
		Rappelons que $\cosh\dist(\isome{f}{}{}(\ell),\ell)=\isome{f}{}{}(\ell)\cdot \ell=\deg f$, donc $\deg f\leq \cosh(
		2\dist(w_n,\ell)+\epsilon_n)$. Remarquons que $\dist(w_n,\ell)=\argcosh(w_n\cdot \ell)=\argcosh\sqrt{2}$ (en utilisant (\ref{eq_wn})).
		Finalement, d'après le choix (\ref{hypothese_star}) de $\epsilon_n$, nous majorons le degré de $f$ par $4$ : \[\deg f \leq \cosh(2\argcosh(\sqrt{2})+\epsilon_n)< 4.\]
		
		Il nous reste à montrer que $f$ ne peut pas être de degré $2$ ou $3$. Dans chacun des deux cas nous raisonnons par l'absurde et arrivons à l'inégalité $\dist(\isome{f}{}{}(\ell), w_n)>\argcosh(\sqrt{2})+\epsilon_n$. Or, par l'inégalité triangulaire, nous obtenons : \begin{align*}
				\dist(\isome{f}{}{}(\ell),w_n)&\leq\dist(\isome{f}{}{}(\ell),\isome{f}{}{}(w_n))+\dist(\isome{f}{}{}(w_n),w_n)\\
				&\leq \dist(w_n,\ell)+\epsilon_n = \argcosh(\sqrt{2})+\epsilon_n,
				\end{align*}
		ce qui nous mène à une contradiction. 
		
		Supposons que $\deg(f)=3$, alors $f^{-1}$ possède un point base de multiplicité $2$ et quatre points base de multiplicité $1$, d'après les relations de Noether (voir \cite[p.39]{KS}). Ainsi, l'action de $\isome{f}{}{}$ sur la classe $\ell$ de Picard-Manin est définie par $\isome{f}{}{}(\ell)=3\ell-2e_0-e_1-e_2-e_3-e_4$ avec $e_i$ les classes des diviseurs exceptionnels au-dessus de $\P_{\bar{\K}}$. D'après la remarque \ref{rmq_tousdistincts}, nous obtenons :\[(2e_0+e_1+e_2+e_3+e_4)\cdot r_n\geq \frac{1}{n}(-2(n-1)-(n-1)-1-1-1)=-3,\] ce qui nous donne $-\isome{f}{}{}(\ell)\cdot r_n\geq -3$. Par hypothèse nous avons $\isome{f}{}{}(\ell)\cdot \ell = 3$, ce qui nous permet, en utilisant (\ref{eq_wn}), d'établir l'inégalité : \[\isome{f}{}{}(\ell)\cdot w_n = \sqrt{2}\isome{f}{}{}(\ell)\cdot \ell - \frac{1}{\sqrt{2}}\isome{f}{}{}(\ell)\cdot r_n\geq 3\sqrt{2}-\frac{3}{\sqrt{2}}=\frac{3}{\sqrt{2}}.\] En prenant le cosinus hyperbolique réciproque et en utilisant (\ref{hypothese_star}), nous arrivons à la contradiction attendue : \[\dist(\isome{f}{}{}(\ell),w_n)\geq \argcosh\frac{3}{\sqrt{2}}>\argcosh(\sqrt{2})+\epsilon_n.\]
		
		Supposons maintenant que $\deg f=2$, alors d'après les relations de Noether, $f^{-1}$ possède trois points base de multiplicité $1$. De même que dans le cas précédent, nous considérons l'action de $\isome{f}{}{}$ sur $\ell$; $\isome{f}{}{}(\ell)=2\ell-e_0-e_1-e_2$ avec $e_i$ les classes des diviseurs exceptionnels au-dessus de $\P_{\bar{\K}}$. D'après la remarque \ref{rmq_tousdistincts}, nous avons :\[
		-\isome{f}{}{}(\ell)\cdot r_n = (e_0+e_1+e_2)\cdot r_n\geq \frac{1}{n}(-(n-1)-(n-1) - 1)= -2 +\frac{1}{n}. \] De plus, $\isome{f}{}{}(\ell)\cdot \ell =2$, ce qui implique :\[ \isome{f}{}{}(\ell)\cdot w_n = \sqrt{2}\isome{f}{}{}(\ell)\cdot \ell - \frac{1}{\sqrt{2}}\isome{f}{}{}(\ell)\cdot r_n\geq 2\sqrt{2}-\frac{2}{\sqrt{2}}+\frac{1}{n\sqrt{2}}=\sqrt{2}+\frac{1}{n\sqrt{2}}.\]
		De même que précédemment, nous obtenons la contradiction : \[\dist(\isome{f}{}{}(\ell)\cdot w_n)\geq \argcosh\Big(\sqrt{2}+\frac{1}{n\sqrt{2}}\Big)>\argcosh(\sqrt{2})+\epsilon_n.\qedhere\] 
		\end{proof}

		\begin{proof}[Démonstration du théorème \ref{prop_principale}] Le groupe $\Bir(\P_{\K})$ étant inclus dans le groupe $\Bir(\P_{\bar{\K}})$, il agit sur l'espace $\H_{\bar{\K}}$.
		Nous allons montrer que $\Bir(\P_{\bar{\K}})$ agit discrètement le long de l'axe de $\isome{h}{n}{}$, c'est-à-dire qu'il existe $\epsilon>0$ tel que l'ensemble $\Fix_{\epsilon}\{\isome{h}{n}{2}(w_n),\isome{h}{n}{-2}(w_n)\}\subset \Bir(\P_{\bar{\K}})$ soit fini. Ainsi, $\Fix_{\epsilon}\{\isome{h}{n}{2}(w_n),\isome{h}{n}{-2}(w_n)\}\cap \Bir(\P_{\K})$ sera fini ce qui signifie que le groupe $\Bir(\P_{\K})$ agit discrètement le long de l'axe de $\isome{h}{n}{}$ et le théorème sera démontré.

			Choisissons un $\epsilon_n>0$ vérifiant (\ref{hypothese_star}), comme dans le lemme \ref{lemme_lineaire}. Soit $f$ une application birationnelle de $\P_{\bar{\K}}$ appartenant à $\Fix_{\epsilon_n}\{\isome{h}{n}{2}(w_n),\isome{h}{n}{-2}(w_n)\}$. Par le lemme \ref{lemme_cvxd}.(\ref{lemme_cvxd_fix}), nous avons l'inclusion :
			\[\Fix_{\epsilon_n}\{\isome{h}{n}{2}(w_n),\isome{h}{n}{-2}(w_n)\}\subset \Fix_{\epsilon_n}\{\isome{h}{n}{}(w_n),\isome{h}{n}{-1}(w_n)\}\subset \Fix_{\epsilon_n}\{w_n\}.\]
			
			De ce fait, $f$ appartient à $\Fix_{\epsilon_n}\{w_n\}$ et en utilisant le lemme \ref{lemme_lineaire}, nous obtenons que $f$ est de degré $1$ : \[f \colon [x:y:z]\mapsto[ax+ky+bz:lx+cy+dz:hx+my+jz].\]

		 	Le but est maintenant de trouver des contraintes sur les coefficients de $f$ pour montrer qu'il n'y a qu'un nombre fini de choix possibles pour une telle application. L'application $f$ appartient aussi à $\Fix_{\epsilon_n}\{\isome{h}{n}{}(w_n),\isome{h}{n}{-1}(w_n)\}$ d'où $\dist(\isome{f}{}{}\isome{h}{n}{-1}(w_n),\isome{h}{n}{-1}(w_n))\leq \epsilon_n$. L'application $\isome{h}{n}{-1}$ étant une isométrie nous avons :\[\dist(\isome{h}{n}{}\isome{f}{}{}\isome{h}{n}{-1}(w_n),w_n)\leq \epsilon_n.\]
			De plus, $\isome{h}{n}{}\isome{f}{}{}\isome{h}{n}{-1}=\isome{(h_nfh_n^{-1})}{}{}$ donc d'après le lemme \ref{lemme_lineaire} l'application $h_nfh_n^{-1}$ est un automorphisme. En considérant la courbe $C$ contractée sur $p_0$ par $h_n^{-1}$ (la droite à l'infini) et le fait que $f$ est un automorphisme, nous avons $f(p_0)=p_0$, ce qui implique que $l=h=0$. 
			En utilisant le même argument et le fait que  $\dist(\isome{f}{}{}\isome{h}{n}{}(w_n),\isome{h}{n}{}(w_n))\leq \epsilon_n$ nous obtenons $f(q_0)=q_0$ qui implique $k=m=0$. À présent, nous pouvons supposer $j=1$. Ainsi, $f$ est un automorphisme affine : \[f \colon (x,y)\mapsto (ax+b,cy+d).\] Calculons $h_n f h_n^{-1}$ et $h_n^{-1}f h_n$ et regardons à quelles conditions ce sont des automorphismes de degré $1$ de $\A_{\bar{\K}}^2$. \begin{align*}
				h_n fh_n^{-1} & =(y,y^n-x)\circ(ax+b,cy+d)\circ (x^n-y,x)\\
									 & =(cx+d,(cx+d)^n+ay-ax^n-b)\\
									 & =(cx+d,x^n(c^n-a)+nc^{n-1}dx^{n-1}+\dots+ncd^{n-1}x+ay+d^n-b).
									\end{align*}
			De même, nous avons :
			\[h_n^{-1} f h_n =(y^n(a^n-c)+na^{n-1}by^{n-1}+\dots+nab^{n-1}y+cx+b^n-d,ay+b).\]					
		Comme les applications $h_n f h_n^{-1}$ et $h_n^{-1}f h_n$ sont des automorphismes affines de $\A_{\bar{\K}}^2$, les coefficients $a$ et $c$ de $f$ vérifient les relations : \[c^n=a \text{ et } c=a^n \] Cela signifie que $a$ et $c$ appartiennent à l'ensemble $\mathbb{U}_{n^2-1}$ des racines $(n^2-1)$-ème de l'unité car $f$ est un automorphisme et donc $a$ et $c$ sont non nuls.

		Supposons dans un premier temps que $n\geq 3$. En considérant les coefficients de $x^{n-1}$ et $y^{n-1}$ nous avons : 
		 $nc^{n-1}d=na^{n-1}b=0$. Comme la caractéristique du corps ne divise pas $n$ et que $a$ et $c$ sont non nuls, nous obtenons \[d=b=0.\]
		Pour $n\geq 3$, $f$ est donc de la forme $(x,y)\mapsto (ax,cy)$ avec $a,c\in\mathbb{U}_{n^2-1}$ et $a^n=c$, ainsi l'ensemble $\Fix_{\epsilon_n}\{\isome{h}{n}{2}(w_n),\isome{h}{n}{-2}(w_n)\}$ est fini.

		Si maintenant $n=2$, nous ne pouvons pas conclure directement.
		Par hypothèse, $f$ appartient à $\Fix_{\epsilon_2}\{\isome{h}{2}{2}(w_2),\isome{h}{2}{-2}(w_2)\}$ cela signifie que \[\dist(\isome{f}{}{}\isome{h}{2}{-2}(w_2),\isome{h}{2}{-2}(w_2))\leq \epsilon_2 \text{ et }\dist(\isome{f}{}{}\isome{h}{2}{2}(w_2),\isome{h}{2}{2}(w_2))\leq \epsilon_2.\] Comme $\isome{h}{2}{-2}$ et $\isome{h}{2}{2}$ sont des isométries nous avons en fait :\[\dist(\isome{h}{2}{2}\isome{f}{}{}\isome{h}{2}{-2}(w_2),w_2)\leq \epsilon_2 \text{ et } \dist(\isome{h}{2}{-2}\isome{f}{}{}\isome{h}{2}{2}(w_2),w_2)\leq \epsilon_2.\]
		D'après le lemme \ref{lemme_lineaire}, les applications $h_2^2fh_2^{-2}$ et $h_2^{-2}fh_2^2$ sont de degré $1$. De plus le point $p_0$ (respectivement le point $q_0$) est le seul point base de $h_2^2$ (respectivement de $h_2^{-2}$) appartenant à $\P_{\bar{\K}}$. Par le même argument que précédemment, $h_2 fh_2^{-1}$ préserve $p_0$ et $h_2^{-1} f h_2$ préserve $q_0$. Nous obtenons que $2cd=0$ et $2ab=0$. Comme nous sommes en caractéristique différente de $2$ et que $f$ est un automorphisme nous devons avoir $a$ et $c$ non nuls donc $b=d=0$.
	Finalement, lorsque $n=2$, $f$ est de la forme $(x,y)\mapsto (ax,cy)$ avec $a,c\in\mathbb{U}_{3}$ et $a^2=c$. L'ensemble $\Fix_{\epsilon_2}\{\isome{h}{2}{2}(w_2),\isome{h}{2}{-2}(w_2)\}$ est fini comme attendu.
		\end{proof}

\bibliographystyle{smfalpha}
\bibliography{biblio}
\end{document}